# An Approach to the High-level Maintenance Planning for EMU Trains Based on Simulated Annealing

Boliang Lin [a,1], Ruixi Lin [b]

[a] *School of Traffic and Transportation, Beijing Jiaotong University, Beijing 100044, China*
[b] *Department of Electrical Engineering, Stanford University, Stanford, CA 94305, USA*

**Abstract:** A high-speed train needs high-level maintenance when its accumulated running mileage or time reaches predefined threshold. The date of delivering an Electric Multiple Unit (EMU) train to maintenance ranges within a time window rather than be a fixed date. Obviously, changing the delivering date always means a different impact on the supply of EMU trains and operation cost. Therefore, the delivering plan has the potential to be optimized. This paper formulates the EMU train high-level maintenance planning problem as a non-linear 0-1 programming model. The model aims at minimizing the mileage loss of all EMU trains with the consideration of the maintenance capacity of the workshop and maintenance ratio at different times. The number of trains under maintenance not only depends on the current maintenance plan, but also influenced by the trains whose maintenance time span from the last planning horizon to current horizon. A state function is established to describe whether a train is under maintenance. By using this function the constraint of restricting the total number of trains that are under maintenance can be formulated reasonably well. Finally, a simulated annealing algorithm is proposed for solving the problem.

## 1. Introduction

High-level maintenance planning (HMP) for Electric Multiple Unit (EMU) trains is one of the most studied problems and represents an important task in the high-speed railway operation and management.

With the speed up of industrialization process, the air pollution is more and more serious. Among the many sources of air pollution, the transportation industry is a large one. In the mid-long distance transportation, compared to highway and air transport, high-speed railway are better choices for developing green transportation because of its low emissions. In this situation, more and more countries began to pay much attention to the development of high-speed rail.

Take China as example, the high-speed railway mileage had reached more than 22,000 kilometers[2] by the end of 2016, which is the longest in the world. At present, there are 4,665 EMU trains[3] in operation per day, serving about four million passengers in average. In order to complete the task, there about 2,600 EMU trains in the fleet. The maintenance cost of an EMU train with 8 cars is about 15 million CNY per year[4], about one-tenth of the EMU train's price. Therefore, the total cost of maintaining all the EMU

---

trains is about 39 billion CNY in one year.

It is well known that high-speed trains should be inspected/maintained after a given travel distance or time for safety reasons. According to the EMU inspection regulation of CRH series in China, the EMU train maintenance is divided into five levels, of which the first and the second levels are called operational maintenance, while the others are called high-level maintenance.

With the expansion of high-speed railway network, more and more trains are introduced to the fleet. As their accumulated running mileage is growing over time, an increasing number of EMU trains are reaching their maintenance cycles, indicating that the workload at the maintenance facility gets significantly heavier. Thus, it has become an urgent scientific problem to design effective and efficient maintenance plans for the EMU trains.

The desired number of EMU trains to execute passenger transport task is fluctuant during different time periods in a year for a high-speed rail system. For example, during the Spring Festival transport rush which usually lasts for forty days, many additional EMU trains are put into operation, and the trains are rarely allowed to carry out maintenance in this situation. To meet passenger traffic peak demand, EMU trains are usually bring forward to perform high-level maintenance, though this will result in a waste of remaining mileage in terms of the target value.

Therefore, the purpose of optimizing the HMP is to reduce the maintenance frequency reasonably while ensuring safe operations. To this end, the maintenance delivering date should never deviate from the target date too far. Because each EMU train has a time-window during which any day it can be sent to the workshop, in a certain day, the number of trains under maintenance is limited due to the maintenance rate which is various in different period. Therefore, in mathematics, the high-level maintenance plan is a combinatorial optimization problem. Hence, HMP is a scientific problem that needs to be addressed urgently both in theory and practice.

The remainder of this paper is organized as follows. In Section 2, a literature review of recent studies on train maintenance planning related problems are presented. Section 3 gives more details and a formal definition for the HMP problem. A 0-1 nonlinear programming model is proposed to mathematically describe the HMP problem in Section 4. Section 5 develops a simulated annealing based solution framework. Finally, conclusions and research prospects are drawn in Section 6.

## 2. Literature Review

For the optimization of train maintenance plan, the existing literatures are more concerned about the operational-level maintenance schedule. Maróti et al. (2005, 2007) developed a multi-commodity flow type model for preventive maintenance routing of train units and an integer programming model for the urgent preventive maintenance. Wang et al. (2010) built an integer linear programming model with maximum traveling mileage of train units. Wang et al. (2012) established an integer programming model for optimizing EMU train assignment and maintenance schedule. Lin et al. (2013) described three different operation modes of EMU train based on the characteristics of operational maintenance on fixed and unfixed train routing. Giacco et al. (2014) proposed a two-step approach and a mixed-integer linear-programming formulation for integrating short-term maintenance planning respectively. Lai et al. (2015) developed an exact optimization model to improve the efficiency in rolling stock usage with consideration of all necessary regulations and practical constraints. And a hybrid heuristic process was also developed to improve solution quality and efficiency. Empirical results demonstrated that the heuristic process can successfully increase the efficiency of rolling stock use by about 5%. Li et al. (2016) proposed a 0-1 integer programming model for EMU assignment and maintenance scheduling. These literatures seemly did not include the high-level maintenance schedule of EMU trains.

Similar to EMU maintenance scheduling, scheduled maintenance planning problems have been studied in other fields. Moudani et al. (2000) proposed a comprehensive process mixed a dynamic programming approach and a heuristic technique to solve the joint problem of fleet allocation and maintenance scheduling. Sriram et al. (2003) presented an innovative formulation for the aircraft re-assignment and maintenance scheduling and a heuristic method to solve the problem efficiently and quickly. Keysan et al. (2010) addressed both tactical and operational planning for scheduled maintenance of air transportation. Deris et al. (1999) modelled ship maintenance scheduling as a constraint satisfaction problem. Haghani et al. (2002) dealt with the problem of scheduling bus maintenance activities. Kralj et al. (1995) described a multi-objective optimization approach for annual preventive maintenance scheduling of fossil fuel thermal units in a power system. Dahal et al. (2007) presented a comparative study of the applications of a number of GA-based approaches to solve the generator maintenance scheduling problem in power systems using a reliability criterion.

In addition, maintenance schedule is also a periodic scheduling problem, because maintenance activities have to be performed at regular intervals. Periodic scheduling is a well-researched area with a broad range of applications. Grigoriev et al. (2006) introduced an integer linear program formulation for scheduling periodic maintenance. Chen (2009) minimized the number of tardy jobs subject to periodic maintenance and non-resemble jobs. Benmansoura et al. (2014) minimized the weighted sum of maximum earliness and maximum tardiness costs. Phan et al. (2015) developed a multi-stage optimization framework for determining periodic inspection intervals for geo-distributed infrastructure systems subject to hidden failures. Moreover, few studies have also researched the maintenance interval for reliability, such as Sriskandarajah et al. (1998) and Moghaddam et al. (2011).

In the field of EMU train high-level maintenance research, Li et al. (2013) developed a method to forecast the maintenance quantity of EMU trains in arbitrary time period accurately and efficiently. Wang et al. (2016) established two integer programming models for the Repair Shop Scheduling Problem of EMU. To the best of our knowledge, very few studies address the high-level maintenance plan for EMU trains.

## 3. Problem Description

For safety reasons, high-speed trains must be inspected after a given travel distance or time. There are five levels of inspection/maintenance which are labeled as level Ⅰ ~ Ⅴ. Take CRH2 series train as example, the inspection regulation is shown in Table 1.

Table 1. Inspection regulations for CRH2 series train

| Level | Interval requirements | Maintenance content |
|-------|----------------------|---------------------|
| Ⅰ | 4000±400 km or 48 hours | Daily inspection, consumables replacement, and cleaning operation |
| Ⅱ | $30^{+3}_{-3}$ thousand km or one month | Monthly inspection, including hollow shaft ultrasonic flaw detection, gearbox change oil and wheel tread repair. |
| Ⅲ | $600^{+20}_{-50}$ thousand km or 1.5 years | Bogie inspection and maintenance |
| Ⅳ | $1200^{+50}_{-100}$ thousand km or 3 years | To disassemble the each subsystem of the EMU for inspection and maintenance and the EMU is paint again. |
| Ⅴ | $2400^{+100}_{-100}$ thousand km or 6 years | The EMU needs overhauling and most of EMU parts are replaced. |

The level I and II are both operational inspection with short cycle. The remainder belong to high-level maintenance with long cycle. In general, the operational maintenance plan is co-optimization with the rolling stock assignment plan at a tactical level, and there is already a lot of literature for this field. High-level maintenance is usually of long interval and the maintenance time is also relatively long. Hence, HMP will influence the demand of using EMU in some periods such as Spring Festival travel rush and MMS should be optimized from the strategic level annually or with longer planning horizon.

The interval of mileage and time between high-level maintenance of adjacent levels of CRH2 series train are shown in Fig. 1, where yellow triangle, pink hexagons and five red pentagram denote level III, IV and V, respectively. The symbol "600/1.5" means 600 thousand kilometers or 1.5 years. When a new EMU has run for 600 thousand km or in operation for 1.5 years, it will meet its first inspection at level III. After another 600 thousand km, that is to say, after a total of 1.2 million km or three years, the train will need its first inspection at level IV. After a total of 1.8 million km or 4.5 years, the EMU will meet second maintenance of level III, and so on.

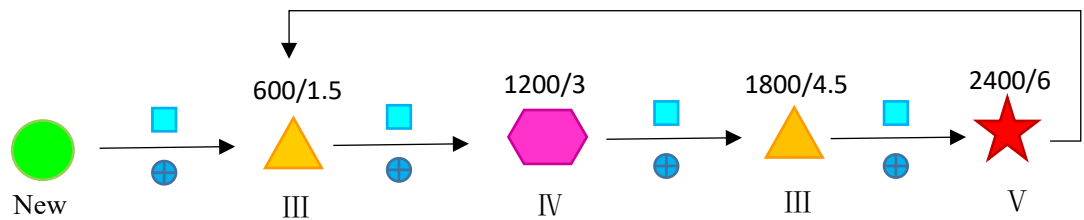

Fig. 1. High-level maintenance intervals for the CRH2 series train

It's easy to find that the interval between two adjacent levels of high-level maintenance is 600 thousand km or 1.5 years for CRH2 series train. Besides, level I (represented by dark blue ⊕ symbol) and level II (represented by light blue square) are implemented between two adjacent levels of high-level maintenance.

It is necessary to explain that the maintenance procedures of higher-level cover all contents of lower-level maintenance. Therefore, after one level of maintenance process, all accumulated running mileage and days associated with that level and the corresponding lower-level classes of maintenance are set to zero.

According to the actual operating performance, the annual running mileage of CRH2 series EMU train is more than 500,000 km. Therefore, their cumulative running mileage usual reaches the target number earlier than operating days. Based on such a fact, we could take into account the cumulative mileage only and ignore the time interval in this paper. In the EMU management information system, running mileage and maintenance track and other data for all trains are well documented. We can easily deduce the inspection expired date $t_m^{\text{expired}}$ and maintenance level $g_m$ for EMU train $m$ on these data. According to inspection regulations, there is an allowance on both sides of the threshold. In other words, lower bound of the window is earlier than $t_m^{\text{expired}}$, and the upper bound is later than $t_m^{\text{expired}}$. Take level III of CRH2 series train as example, the mileage interval is $600^{+20}_{-50}$, which means that the fluctuation range of cumulative operating mileage can be described as [550,000 km, 620,000 km]. Therefore, each train has a time-window during which any day it can be sent to maintain, as shown in Fig. 2.

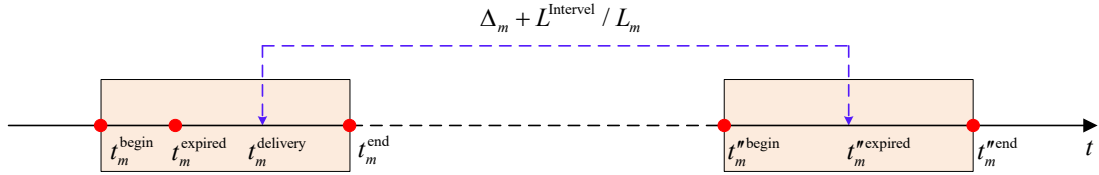

Fig. 2. The time-window during which any day train can be sent to maintain

It is obvious that the starting time and the width of a time-window might be different for each EMU train. Generally, a higher maintenance level represents a longer time window. The upper bound (ending time) $t_m^{\text{end}}$ and the lower bound (starting time) $t_m^{\text{begin}}$ of a time window can be calculated as following equations:

$$t_m^{\text{begin}} = t_m^{\text{expired}} - L_{g_m}^{\text{left}} / L_m \qquad (1)$$

$$t_m^{\text{end}} = t_m^{\text{expired}} + L_{g_m}^{\text{right}} / L_m \qquad (2)$$

where $L_{g_m}^{\text{left}}$ and $L_{g_m}^{\text{right}}$ denote the negative and positive offset value of the target mileage distance from the next maintenance, $g_m$ is maintenance level, and $L_m$ denotes the daily running mileage for train $m$. If an EMU train is delivered for maintenance at the early stage of planning horizon, there may be another high-level maintenance arranged for this train during the same horizon. In other words, the second delivery time window might occur in the planning horizon (time window on the right-hand side in Fig. 2). Let $t_m^{\text{delivery}}$ be the deliver time, then the theoretical expired date of next maintenance will be $t_m''^{\text{expired}} = t_m^{\text{delivery}} + \Delta_m + L^{\text{interval}} / L_m$, where $L^{\text{interval}}$ (e.g. 600,000 km) represents the interval between two adjacent levels of high-level maintenance.

In order to better understand the computing of $t_m^{\text{begin}}$ and $t_m^{\text{end}}$, we would like to give more explanations of the time window using a three-trains example (see Fig. 3). For the first train (number EMU_001), we assume $t_1^{\text{expired}} = 127$, $g_1 = \text{III}$, $L_{g_1}^{\text{left}} = 50,000$ km, $L_{g_1}^{\text{right}}$ =20,000 km, $L_1$ =1,600 km, we then have $t_1^{\text{begin}} = 127 - 50/1.6 = 95.75$ and $t_1^{\text{end}} = 127 + 20/1.6 = 139.5$.

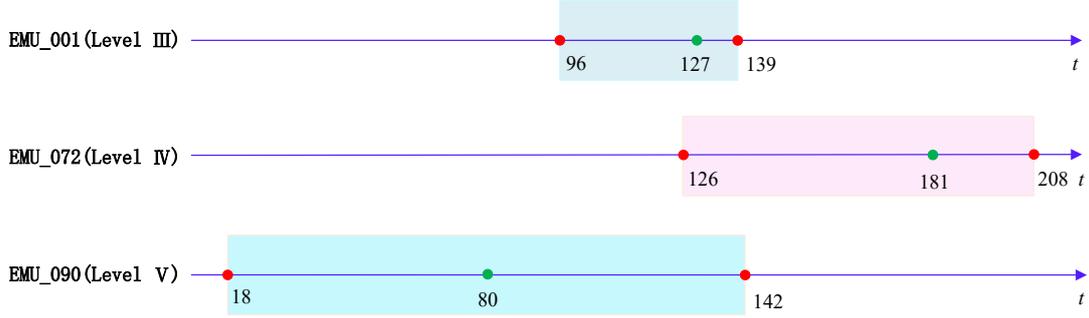

Fig. 3. The lower and upper bound of the time window

The boundary will be rounded toward window center according to the maintenance regulations if the boundary is a fraction. Therefore, the time-window for train EMU_001 is [96, 139] (see Fig. 3). Similarly, for the second train (number EMU_072), we assume $t_{72}^{\text{expired}}$ =181, $g_{72}$ =IV, $L_{g_{72}}^{\text{left}}$ =100,000 km, $L_{g_{72}}^{\text{right}}$ =50,000 km and $L_{72}$ =1,800 km. Hence, $t_{72}^{\text{begin}}$ =181-100/1.8=125.4, $t_{72}^{\text{end}}$ =181+50/1.8=208.7. The time-window for EMU_072 is [126, 208] (see Fig. 3). For the third train (number EMU_090), we assume $t_{90}^{\text{expired}}$ =80, $g_{90}$ =V, $L_{g_{90}}^{\text{left}}$ =100,000 km, $L_{g_{90}}^{\text{right}}$ =100,000 km, $L_{90}$ =1,600 km. Hence, $t_{90}^{\text{begin}}$ =80-100/1.6=17.5, $t_{90}^{\text{end}}$ =80+100/1.6=142.5. The time-window for train EMU_090 is [18, 142] (see Fig. 3).

The span of planning horizon can be calculated by the following formula:

$$T = \max_m \left\{ t_m^{\text{end}} \right\} \tag{3}$$

For a high-speed railway system, if there are enough available EMU trains and we do not consider the maintenance rate, the ideal delivery date of each train will be the expected expired date. In fact, due to the significantly high purchasing costs of EMU trains, the railway operator generally maintain a reasonable fleet size of EMU trains to meet basic passenger traffic demands. There are two reasons make it impossible to deliver all trains at their corresponding upper bounds of time windows. The first is the non-uniform distribution of expiration date of high-level maintenance. The second is the existence of traffic rush periods, such as the Spring Festival and the summer holiday, which results in the unbalanced demand of EMU trains. Thus, how to make the maintenance delivery date close to the upper bound of time window under the condition of satisfying operational requirements in various time segments, is a problem with great value in practice. Theoretically, the delivery time can be regarded as discrete point in time window (denoted as a binary decision variable). And time windows of different trains often partially overlap with each other, hence, the start time of inspection for different EMU train is interrelated, especially under the constraints of a certain maintenance rate and capacity. Therefore, HMS is mathematically a combinatorial optimization problem with nonlinear constraints.

# 4. Model Formulation for the HMP

The proposed HMS formulation aims to assign a maintenance schedule to each EMU train while maximizing the utilization of remaining mileage. The date for delivering a train to the workshop should be chosen among the discrete nodes within the time window. Also, the current number of trains in maintenance should meet the limitation of maintenance rate and the workshop's capacity. Because the HMP involves a long planning horizon, usual a year, some trains may last from the last planning horizon to the current horizon, we need to deduct the capacity occupied by them or preload them onto the current horizon. Besides, some trains may occasionally be sent to maintain twice a year. Hence, their second time window should be considered.

## 4.1. List of Notations

The parameters, sets and decision variables used in the optimization model are listed in Table 2.

Table 2. Notations used in the formulation

| Symbol | Definition |
|---|---|
| **Parameters** | |
| $t_m^{\text{expired}}$ | The date when the EMU train $m$ should be arranged for maintenance in theory |
| $g_m$ | The maintenance level planned for EMU train $m$ |
| $N_m^{\text{size}}$ | The number of units grouped in EMU train $m$, a train unit composes of eight passenger cars |
| $L_m$ | The daily running mileage of EMU train $m$ |
| $C^{\text{maintain}}$ | The apportioned high-level maintenance cost per train-mileage |
| $C^{\text{income}}$ | The average passenger revenue per train-mileage |
| $\gamma^{\text{SpringRush}}$ | The allowable maintenance rate during the Spring Festival rush, which is equal to the rate of the number of EMU trains in high-level maintenance to the total number of EMU trains |
| $\gamma^{\text{SummerRush}}$ | The allowable maintenance rate during the summer holiday rush |
| $\gamma_i$ | The allowable maintenance rate in the $i$ th time period during the planning horizon (e.g. $i$=2 means the Spring Festival rush and $i$=4 represents the summer holiday rush) |
| $\gamma$ | The maintenance rate in the non-special time period (i.e. outside the rush periods) |
| $\lambda$ | The profit rate of passenger ticket revenue |
| $N^{\text{accept}}$ | The acceptable EMU trains for a workshop within a single day |
| $\Delta_m$ | The maintenance service time of the level $g_m$ for EMU train $m$ |
| $T_b^{\text{SpringRush}}$ | The begin time of the Spring Festival travel rush |
| $T_e^{\text{SpringRush}}$ | The end time of the Spring Festival travel rush |
| $T_b^{\text{SummerRush}}$ | The begin time of the summer holiday travel rush |
| $T_e^{\text{SummerRush}}$ | The end time of the summer holiday travel rush |
| $T$ | The length of the planning horizon |
| $C_k^{\text{capacity}}$ | Maintenance capacity for level $k$ |
| **Sets** | |
| $S^{\text{emu}}$ | The set of all EMU trains |
| $S_m^{\text{window}}$ | The time window of EMU train $m$, during which the EMU train can be delivered to maintain |
| $S^{\text{level}}$ | The set of maintenance levels, which is described as $S^{\text{level}} = \{3,4,5\}$ in this paper |
| $T_i$ | The $i^{\text{th}}$ time period in the planning horizon, which is described as $T_i \subset [0,T]$ here |
| **Decision Variables** | |
| $x_m^t$ | Binary decision variables. It takes the value of one if the EMU train $m$ is delivered for maintenance at the $t^{\text{th}}$ day, and is zero otherwise. |

## 4.2. The State Function

The EMU trains that involve maintenance work within the planning horizon [0, T] consist of three parts: (1) both the maintenance start time and end time fall within the planning horizon; (2) the maintenance start time is before the planning horizon and the end time falls within the planning horizon and (3) the maintenance start time falls within the planning horizon and the end time is after the planning horizon (see Fig. 4).

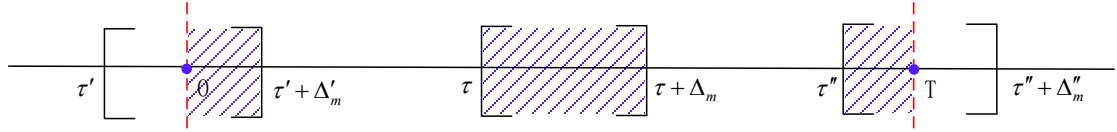

Fig. 4. Time intervals indicating the train is under maintenance

Consider $x_m^\tau = 1$ and the maintenance time interval is $[\tau, \tau + \Delta_m]$ for train $m$. If the start time of last maintenance $\tau' \leq 0$ yet the train's maintenance end time $\tau' + \Delta_m' > 0$, then the time period $[0, \tau' + \Delta_m']$ still impacts on the maintenance rate of the current planning horizon, where $\Delta_m'$ is the maintenance service time of last high-level maintenance. Similarly, if the start time of last maintenance $\tau' > 0$, then the time period $[\tau', \tau' + \Delta_m']$ will also impact on the maintenance rate of the current planning horizon. Given the maintenance start time $\tau$, the expired time of next maintenance $\tau''$ can be computed by the following equation:

$$\tau'' = \tau + \Delta_m + L^{\text{Interval}} / L_m \quad (4)$$

Since the design of a maintenance plan is based on a rolling horizon framework, the deliver time of next maintenance generally belongs to next planning horizon. However, if the deliver time of current maintenance is at the early stage of the planning horizon, then the expired time of next maintenance probably falls into the current planning horizon. Let $\Delta_m''$ denote the maintenance service time of next maintenance and we define the valid maintenance end time as follows:

$$\tau''' = \min\{T, \tau'' + \Delta_m''\} \tag{5}$$

In this way, the time interval of next maintenance that impacts on the current plan is $[\tau'', \ \tau''']$. Therefore, the binary state function that describes whether train $m$ on day $t$ is under maintenance can be defined as follows:

$$f_m(t) = \begin{cases} 1 & \text{If } x_m^\tau = 1 \text{ and } t \in [\tau_b', \tau' + \Delta_m'] \bigcup [\tau, \tau + \Delta_m] \bigcup [\tau'', \tau'''], \forall \tau \in S_m^{\text{window}} \\ 0 & \text{Otherwise} \end{cases} \tag{6}$$

where $\tau_b' = \max\{0, \tau'\}$. If $f_m(t)$ takes a value of 1, then the train is under the maintenance state; otherwise, the train is under the operation state.

## 4.3. Mathematical Formulation

Using above-mentioned notation, the HMP problem can now be written with a nonlinear 0-1 programming formulation as follows:

$$\min \sum_{m \in S^{\text{emu}}} \sum_{t \in S_m^{\text{window}}} (t_m^{\text{expired}} - tx_m^t) L_m N_m^{\text{Size}} (C^{\text{maintain}} + \lambda C^{\text{income}})$$

s.t.

$$\sum_{t \in S_m^{\text{window}}} x_m^t = 1, \qquad \forall m \in S^{\text{emu}} \tag{7}$$

$$\frac{1}{|S^{\text{emu}}|} \sum_{m \in S^{\text{EMU}}} f_m(t) \le \gamma^{\text{SpringRush}}, \quad \forall t \in T^{\text{SpringRush}} \tag{8}$$

$$\frac{1}{|S^{\text{emu}}|} \sum_{m \in S^{\text{emu}}} f_m(t) \le \gamma^{\text{SummerRush}}, \quad \forall t \in T^{\text{SummerRush}} \tag{9}$$

$$\frac{1}{|S^{\text{emu}}|} \sum_{m \in S^{\text{emu}}} f_m(t) \le \gamma, \quad \forall t \in [0, T] \setminus (T^{\text{SpringRush}} \bigcup T^{\text{SummerRush}}) \tag{10}$$

$$\sum_{m \in S^{\text{emu}}} x_m^t \le N^{\text{accept}}, \ \ \forall t \in [0, T] \tag{11}$$

$$\sum_{m \in S^{\text{emu}}} f_m(t) \delta_k^{g_m} \le C_k^{\text{capacity}}, \ \ \forall t \in [0, T], \ \forall k \in S^{\text{level}} \tag{12}$$

$$x_m^t \in \{0, 1\}, \qquad \forall m \in S^{\text{emu}}, \ \ \forall t \in [0, T] \tag{13}$$

In the model, the objective is to minimize the total costs caused by having not sufficiently utilized the remaining mileage. If a train with positive remaining mileage undergoes a maintenance operation, then a portion of transport capacity is wasted and moreover some spare parts are replaced before fully utilizing their useful lives.

Constrains (7) ensure that a train should be sent to maintain exactly once in the related time window. Constrains (8) ~ (10) guarantee that the maintenance rate should be less than the predefined values during the spring peak, summer peak and normal period, respectively. Constrains (11) are to avoid delivering trains into maintenance work too intensively. For example, due to limited maintenance resources in some workshops, only one train can be accepted to proceed maintenance work within a single day. While constraints (12) are workshop capacity restrictions in terms of maintenance levels, i.e., for each maintenance level, the total number of trains under maintenance should be less than a given upper bound. Constraints (13) are the binary constraints on variables.

In constraints (12), the $\delta$ function is defined as follows:

$$\delta_i^j = \begin{cases} 1 & \text{if } i=j \\ 0 & \text{Otherwise} \end{cases} \tag{14}$$

And time periods $T^{\text{SpringRush}}$ and $T^{\text{SummerRush}}$ are the Spring Festival transport rush and the summer holiday transport rush respectively, which are defined as follows:

$$T^{\text{SpringRush}} = [T_b^{\text{SpringRush}}, T_e^{\text{SpringRush}}] \tag{15}$$

$$T^{\text{SummerRush}} = [T_b^{\text{SummerRush}}, T_e^{\text{SummerRush}}] \tag{16}$$

The Spring Festival transport rush usually lasts forty days, many additional high-speed trains are put into operation provisionally, and EMU trains are rarely allowed to go for maintenance in situation.

In fact, constraint (8), (9) and (10) in the above model can be described by the following general formulations:

$$\frac{1}{|S^{\text{emu}}|} \sum_{m \in S^{\text{emu}}} f_m(t) \le \gamma_i, \qquad \forall t \in T_i, \ \ i=1,2,\cdots,n \tag{17}$$

where

$$[0,T] = T_1 \bigcup T_2 \bigcup \cdots \bigcup T_i \cdots \bigcup T_n \tag{18}$$

If the workshop maintenance capacity is independent of the maintenance level, constraints (12) can be simplified as follows:

$$\sum_{m \in S^{\text{emu}}} f_m(t) \le C^{\text{capacity}}, \ \ \forall t \in [0,T] \tag{19}$$

where $C^{\text{capacity}}$ is the total high-level maintenance capacity.

## 4.4. Complexity of the HMP model

Let $m$ denote the fleet size of EMU trains in a high-speed railway system, $n$ be the average width of time windows. Then the total number of 0-1 decision variables is $m \times n$ with the computational complexity of $O(2^{m \times n})$. The number of constraints is $|S^{\text{emu}}|$ for Constrains (7), $T$ for Constraints (8) ~ (10), $T$ for Constraints (11), and $|S^{\text{level}}| \times T = 3T$ for Constraints (12). Therefore, the total number of constraints is $|S^{\text{emu}}| + 5T$. If we do not consider the differences among the maintenance levels, then the number of constraints for workshop capacity constrains is $T$ and the total number of constraints will be $|S^{\text{emu}}| + 3T$.

# 5. Simulated Annealing Based Solution Approach

Simulated annealing (SA) method is employed for solving the model. The SA method was proposed independently by Kirkpatrick et al. (1983) and Černy (1985). The motivation for the SA method was to solve combinatorial optimization problems, and the details of SA method can be referred to Kirkpatrick *et al.* (1983). The steps for applying the SA method in our model solution are as follows (see Lin et al. 2012).

## 5.1. Energy Function

Note that the average cost per train-km is independent of train individuals, the cost parameter $C^{\text{maintain}} + \lambda C^{\text{income}}$ is thus a constant in the objective function. Therefore, we can re-define the minimal function of remaining mileage as follows:

$$H^0(X) = \sum_{m \in S^{\text{emu}}} \sum_{t \in S_m^{\text{window}}} (t_m^{\text{expired}} - tx_m^t) L_m N_m^{\text{size}} \tag{20}$$

where $X = (\cdots, x_m^t, \cdots)$ represents the solution vector.

Because constraints associated with the HMP model are expressed as equalities and inequalities, it is difficult to obtain a feasible solution for a real-world railway operator with a considerably large fleet size. However, in the HMP problem, the constraints can be classified into two groups; "easy" constraints and "difficult" constraints. Obviously, constraints (7) are easy constraints and (11), (17) and (19) belong to the difficult class. Thus, we can convert the difficult constraints into a sequence of penalty functions.

The maintenance rate constraint is converted into the following penalty function:

$$H^1(X) = \max \left\{ 0, \frac{1}{|S^{\text{emu}}|} \sum_{m \in S^{\text{emu}}} f_m(t) - \gamma_i \right\} \tag{21}$$

The penalty function for the daily acceptable number of EMU trains constraint is written by:

$$H^2(X) = \max \left\{ 0, \sum_{m \in S^{\text{emu}}} x_m^t - N^{\text{accept}} \right\} \tag{22}$$

While the workshop capacity constraint penalty function can be expressed as follows:

$$H^3(X) = \max \left\{ 0, \sum_{m \in S^{\text{emu}}} f_m(t) - C^{\text{capacity}} \right\} \tag{23}$$

Therefore, the energy function $E(X)$ of the HMP model can be defined as follows:

$$E(X) = H^0(X) + \beta_1 \sum_{i=1}^{n} \sum_{t \in T_i} H^1(X) + \beta_2 \sum_{t \in [0,T]} H^2(X) + \beta_3 \sum_{t \in [0,T]} H^3(X) \tag{24}$$

In this way, the HMP model can be rewritten in the following form that is suitable for applying the simulated annealing structure.

$$\min \quad E(X)$$

s.t.

$$\sum_{t \in S_m^{\text{window}}} x_m^t = 1, \qquad \forall m \in S^{\text{emu}} \tag{25}$$

$$x_m^t \in \{0,1\}, \qquad \forall m \in S^{\text{emu}}, \quad \forall t \in [0,T] \tag{26}$$

## 5.2. Solution of the HMP

To use the simulated annealing method for the HMP problem, one needs a set of possible solutions. Let $\Omega = \{ X \mid \text{meets (25) and (26)} \}$ denote the set of solutions of the model. The minimum of $E(X)$ over the set $\Omega$ is sought if an optimal solution to the HMP problem is obtained once an element, called a global minimum, $X^* \in \Omega$ with the property $E(X^*) \leq E(X)$ for all $X \in \Omega$ is found.

## 5.3. Initial Solution

An initial solution $X_0 = (\cdots, x_m^t, \cdots)$ can be generated by the following two methods. The easiest way to generate an initial solution is to deliver the trains to maintenance on their theoretical expired dates, i.e. $x_m^t = t_m^{\text{expired}}$. Another method is to simply adopt the manual solution designed by highly experienced dispatchers as the initial solution. In theory, both initial solutions hardly impact on the optimality of final solutions.

## 5.4. Neighborhood Solution

Noted that the current solution of the HMP model is $X_j$, we can select one train $m \in S^{\text{emu}}$ randomly; if $x_m^t = 1$, we change its value to zero. Then we select another date $t' \in S_m^{\text{window}}$ from the time window and set $x_m^{t'} = 1$. In this way, a new HMP solution $X_{j+1}$ is created.

## 5.5. Simulated Annealing Structure

The strategy implemented by SA consists of exploring the solution space starting from an arbitrary selected solution, and then generating a new one by perturbing it. Every time a new solution is generated, then its cost is evaluated and the new solution is either accepted or rejected according to an acceptance rule. The general algorithm of SA can be described by the following flow chart.

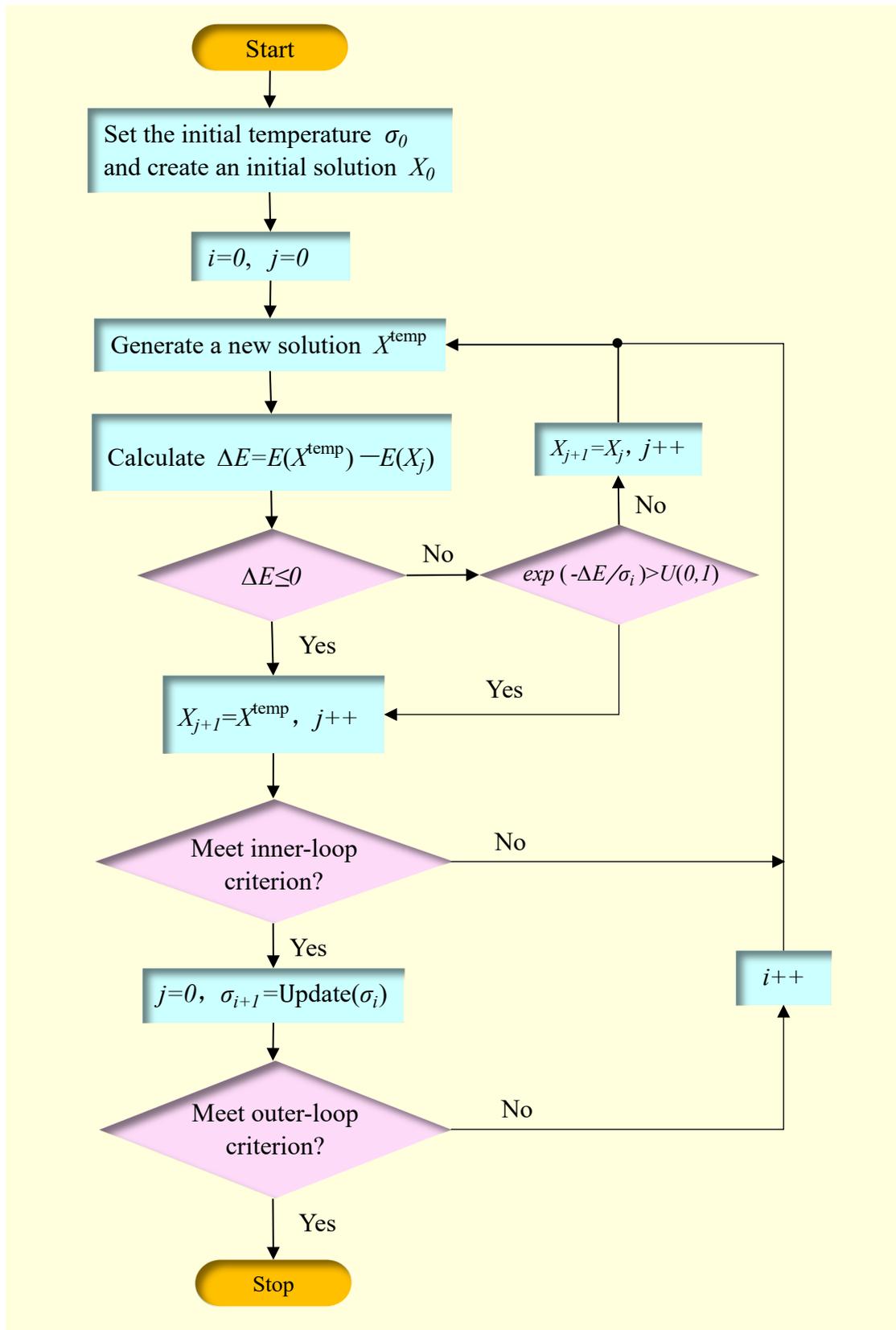

Fig. 5. Flow chart of the simulated annealing based solution approach.

The SA algorithm above needs to specify the inner loop, outer loop criteria, temperature update function and initial temperature.

- Inner-loop criterion (equilibrium condition): Let $N_i^{\text{generated}}$ be the number of new solutions generated, $N_i^{\text{accepted}}$ be the number of accepted solutions at temperature $\sigma_i$. When the condition $N_i^{\text{generated}} \geq h^1 |X|$ or $N_i^{\text{accepted}} \geq h^2 |X|$ is met, the inner-loop criterion is met. In this paper, we set the coefficients $h^1 = 3$ and $h^2 = 6$.
- Outer-loop criterion (convergence criterion): Execution of the algorithm is stopped when the acceptance rate is less than a threshold (e.g., $\varepsilon = 0.001$) or the average cost does not change significantly for consecutive values of temperature. For example, SA is stopped when the value of temperature function does not change significantly for 30 times of consecutive cooling.
- Temperature update function (decrement rule): The update function stands for the decrement rule of the control temperature. From Kirkpatrick's annealing perspective (1983), the temperature updated function has the geometric rate:

$$\sigma_i = h^3 \sigma_{i+1} \qquad (27)$$

  The temperature update coefficient $h^3$ is called the cooling rate and is set as 0.97 in this paper.
- Initial temperature: The value of initial temperature is chosen so that the corresponding acceptance probability density is relatively close to the density for $T = \infty$.

## 6. Conclusions

This paper investigates the high-level maintenance planning problem for EMU trains emerged in the field of high-speed railway operation and management. A 0-1 non-linear programming model is proposed for the problem. The objective of the model aims to minimize the total costs caused by having not sufficiently utilized the remaining mileage, and the constrains include unique deliver date requirements, maintenance rate restrictions and workshop capacity. We consider various time periods that have different maintenance rate requirements in a year, such as the Spring Festival rush and the summer holiday rush. To describe whether an EMU train is under maintenance, a state function is designed. With the help of the function, the maintenance rate constraints can be nicely formulated. Because of the NP-complete nature of the HMP problem, we focus on a simulated annealing based heuristics rather than exact solution algorithms to solve the problem. We divide the constraints into easy ones and difficult ones and convert the later into a penalty function. Furthermore, we outline the general framework of the solution procedure and some details are also provided. Our future work is to test and apply the proposed approach to real-world large-scale problem instances.